\documentclass[preprint,11pt]{elsarticle}

\usepackage{amsmath,amsthm,amsfonts,amssymb,latexsym,mathrsfs,color}
\usepackage{hyperref}
\usepackage{tikz}
\usepackage{bm}

\newtheorem{theorem}{Theorem}%[section]

\newtheorem{example}[theorem]{Example}

\newcommand{\lrf}[1]{\lfloor #1\rfloor}

\journal{}

\begin{document}
\begin{frontmatter}



\title{An interlacing pattern between the types $C_n$ and $D_n$ coordinator polynomials}

\author[focal]{Jun-Ying Liu}
\ead{jyliu6@163.com}
\author[focal]{Shi-Mei Ma}
\ead{shimeimapapers@163.com}
%\cortext[cor1]{Corresponding author: Shi-Mei Ma}
\address[focal]{School of Mathematics and Statistics, Shandong University of Technology, Zibo, Shandong 255000, P.R. China}
\begin{abstract}
The study of coordinator polynomials of Weyl group lattices was initiated by Conway and Sloane.
In 2013, by using a trigonometric substitution approach, Wang and Zhao proved the real-rootedness of the type $D$ coordinator polynomials. Subsequently, 
Xie and Zhang proved the compatibility between the types $C_n$ and $D_n$ coordinator polynomials. 
In this paper, we discover the interlacing pattern between the types $C_n$ and $D_n$ coordinator polynomials, which may be named as the second pattern of interlacing zeros.
\end{abstract}

% %A polynomial $f(x)$ is called bi-gamma-positive if it can be decomposed as $f(x)=a(x)+xb(x)$, 
%%where $a(x)$ and $b(x)$ are both gamma-positive satisfying $\deg a(x)=1+\deg b(x)$ or $b(x)=0$. 
%%A polynomial is said to be ratio monotone if its coefficient sequence is ratio monotone.
% which settles a strong version of an open problem of Brenti in 1994

\end{frontmatter}

\section{Introduction}
%%%%%%%%%%%%%%%%%%%%%%%%%%%%%%%%%%%%%%%%%%%
%%%%%%%%%%%%%%%%%%%%%%%%%%%%%%%%%%%%%%%%%%%%%%%%%%%%%%%%%%%%%%%%%%%%%%%%%%%%%%%%%
%%%%%%%%%%%%%%%%%%%%%%%%%%%%%%%%%%%%%%%%%%%%%%%%%%%%%%%%%%%%%%%%%%%%%%%%%%%%%%%%%
%%%%%%%%%%%%%%%%%%%%%%%%%%%%%%%%%%%%%%%%%%%%%%%%%%%%%%%%%%%%%%%%%%%%%%%%%%%%5%%%%
Let $\mathbf{e}_i$ denote a vector with the $i$th entry one and all other entries zeros.
The root lattices $A_n,~B_n,~C_n$ and $D_n$ can be defined to be generated as
monoids respectively by
\begin{align*}
M_{A_n}&=\left\{\pm(\mathbf{e}_i-\mathbf{e}_j)| 1\leqslant i<j\leqslant n+1 \right\},\\
M_{B_n}&=\left\{\pm\mathbf{e}_i\pm\mathbf{e}_j| 1\leqslant i<j\leqslant n \right\}\cup \left\{\pm\mathbf{e}_i| 1\leqslant i\leqslant n \right\},\\
M_{C_n}&=\left\{\pm\mathbf{e}_i\pm\mathbf{e}_j| 1\leqslant i<j\leqslant n \right\}\cup \left\{\pm2\mathbf{e}_i| 1\leqslant i\leqslant n \right\},\\
M_{D_n}&=\left\{\pm\mathbf{e}_i\pm\mathbf{e}_j| 1\leqslant i<j\leqslant n \right\}.
\end{align*}

The study of coordinator polynomials of rooted lattices was initiated by Conway and Sloane~\cite{Conway97}. They found that 
the coordinator polynomials of type $A_n$ are given as follows:
$$h_{A_n}(x)=\sum_{k=0}^n\binom{n}{k}^2x^k,$$
which is now known as the type $B$ Narayana poolynomials.
As pointed out by Conway and Sloane~\cite{Conway97}, $$h_{A_n}(x)=(1-x)^nL_n\left(\frac{1+x}{1-x}\right),$$
where $L_n(x)$ is the $n$th Legendre polynomial, and so $h_{A_n}(x)$ is real-rootedness.
Subsequently, the explicit expressions for coordinator polynomials of type $B_n$, $C_n$ and $D_n$ are discovered:
 \begin{align*}
h_{B_n}(x)&=\sum_{k=0}^n\binom{2n+1}{2k}x^k-2nx(1+x)^{n-1};\\
h_{C_n}(x)&=\sum_{k=0}^n\binom{2n}{2k}x^k;\\
h_{D_n}(x)&=\frac{(1+\sqrt{x})^{2n}+(1-\sqrt{x})^{2n}}{2}-2nx(1+x)^{n-2}.
\end{align*}
It should be noted that
\begin{equation}\label{hD}
h_{D_n}(x)=h_{C_n}(x)-2nx(1+x)^{n-2}.
\end{equation}
According to~\cite[A086645]{Sloane}, one has
\begin{equation*}\label{hCnx-def}
\sum_{k=0}^n\binom{2n}{2k}x^k=\sum_{k=0}^{\lrf{n/2}}4^k\binom{n}{2k}x^k(1+x)^{n-2k}.
\end{equation*}
The Chebyshev polynomials of the first kind can be defined by the identity $T_n(\cos \theta)=\cos(n \theta)$.
Explicitly, $$T_n(x)=\sum_{k=0}^{\lrf{n/2}}\binom{n}{2k}x^{n-2k}(x^2-1)^k,$$
and the zeros of $T_n(x)$ are $x_k=\cos\frac{(2k-1)\pi}{2n}$, where $k=1,2,\ldots,n$, see~\cite{Rivlin90}.
Note that 
$$
h_{C_n}(x)=(1-x)^n T_n\biggl(\frac{1+x}{1-x}\biggr).
$$
Therefore, the zeros of $h_{C_n}(x)$ are given explicitly by
\begin{equation*}\label{zeroC}
x_{n,k}^c=-\tan^2\frac{(2k-1)\pi}{4n},\qquad k=1,2,\ldots,n.
\end{equation*}

Wang and Zhao~\cite{Wang13} found that $h_{B_{16}}(x)$ has 2 non-real roots, and by using a trigonometric substitution approach,
they mainly showed that the zeros of $h_{D_n}(x)$ are given by $x_{n,k}^d=-\tan^2{\frac{\phi_k}{2}}$ 
for some $\phi_k\in \left(\frac{k\pi}{n},\frac{(k+1)\pi}{n}\right)$, $k=0,1,\ldots,n-1$. However, Wang and Zhao's result offers information no 
more than the existence of real zeros of $h_{D_n}(x)$. 
Xie and Zhang~\cite{Xie2014} proved the compatibility between the types $C_n$ and $D_n$ coordinator polynomials.
A natural idea is to explore more information on the distribution of the zeros of $h_{C_n}(x)$ and $h_{D_n}(x)$.

Let $x_{n,1}^c<x_{n,2}^c<\cdots<x_{n,n}^c$ be the zeros of $h_{C_n}(x)$,
and let $x_{n,1}^d<x_{n,2}^d<\cdots<x_{n,n}^d$ be those of $h_{D_n}(x)$.
By~\eqref{hD}, we see that
$h_{C_n}\left(x_{n,k}^d\right)=2nx_{n,k}^d\left(1+x_{n,k}^d\right)^{n-2}$.
Computational evidence supports the following pattern:
$$
x_{n,1}^c<x_{n,1}^d<x_{n,2}^d<x_{n,2}^c<x_{n,3}^c<x_{n,3}^d<x_{n,4}^d<x_{n,4}^c<\cdots,
$$
i.e., non-extreme consecutive pairs of $x_{n,i}^c$ ``interlace'' those of $x_{n,i}^d$, which may be named as the {\it second pattern of interlacing zeros}, since this is different from 
the classical definition of interlacing zeros~\cite{gw96,LWaam}.
We believe that this new interlacing pattern may be used to study similar interlacing properties in other families of real-rooted polynomials.

We can now conclude the main result of this paper.
\begin{theorem}\label{thm:interlacing}
Let $x_{n,1}^c<x_{n,2}^c<\cdots<x_{n,n}^c$ be the zeros of $h_{C_n}(x)$,
and let $x_{n,1}^d<x_{n,2}^d<\cdots<x_{n,n}^d$ be those of $h_{D_n}(x)$.
For $n\geqslant 3$, they satisfy the second pattern of interlacing zeros:
\begin{itemize}
\item [$(i)$] If $n=2m$, then
\begin{align}\label{interlace2m}
x_{n,1}^{c}
<x_{n,1}^{d}
<x_{n,2}^{d}
<x_{n,2}^{c}
<\cdots
<x_{n,n-1}^{c}
<x_{n,n-1}^{d}
<x_{n,n}^{d}
<x_{n,n}^{c};
\end{align}
\item[(ii)] If $n=2m+1$ with even $m$, then
\begin{align}\label{interlace4r+1}
x_{n,1}^{c}
&<x_{n,1}^{d}
<x_{n,2}^{d}
<x_{n,2}^{c}
<\cdots
<x_{n,m}^{d}
<x_{n,m}^{c}
<x_{n,m+1}^{c}
=x_{n,m+1}^{d}
=-1 \notag\\
&<x_{n,m+2}^{c}
<x_{n,m+2}^{d}
<\cdots
<x_{n,n-1}^{c}
<x_{n,n-1}^{d}
<x_{n,n}^{d}
<x_{n,n}^{c};
\end{align}
\item[(iii)]If $n=2m+1$ with odd $m$, then
\begin{align}\label{interlace4r+3}
x_{n,1}^{c}
&<x_{n,1}^{d}
<x_{n,2}^{d}
<x_{n,2}^{c}
<\cdots
<x_{n,m}^{c}
<\bm{x_{n,m}^{d}}
<x_{n,m+1}^{c}
=x_{n,m+1}^{d}
=-1 \notag\\
&<\bm{x_{n,m+2}^{d}}
<x_{n,m+2}^{c}
<\cdots
<x_{n,n-1}^{c}
<x_{n,n-1}^{d}
<x_{n,n}^{d}
<x_{n,n}^{c}.
\end{align}
\end{itemize}
\end{theorem}
%%%%%%%%%%%%%%%%%%%%%%%%%%%%%%%%%%%%%%%%%%%%%%%%%%%%%%%%%%%%%%%%%%%%%%%%%%%%%%%%%
%%%%%%%%%%%%%%%%%%%%%%%%%%%%%%%%%%%%%%%%%%%%%%%%%%%%%%%%%%%%%%%%%%%%%%%%%%%%%%%%%
%%%%%%%%%%%%%%%%%%%%%%%%%%%%%%%%%%%%%%%%%%%
\section{Proof of Theorem~\ref{thm:interlacing}}
%%%%%%%%%%%%%%%%%%%%%%%%%%%%%%%%%%%%%%%%%%%
%%%%%%%%%%%%%%%%%%%%%%%%%%%%%%%%%%%%%%%%%%%%%%%%%%%%%%%%%%%%%%%%%%%%%%%%%%%%%%%%%
%%%%%%%%%%%%%%%%%%%%%%%%%%%%%%%%%%%%%%%%%%%%%%%%%%%%%%%%%%%%%%%%%%%%%%%%%%%%%%%%%
Following the same way as in~\cite{Wang13}, 
let
$x=-y^2$ and $y=\tan\frac{\phi}{2}$, where $\phi\in(0,\pi)$. 
By the Weierstrass substitution, we have
$$\sin\phi=\frac{2y}{1+y^2},~\cos\phi=\frac{1-y^2}{1+y^2}.$$
Wang and Zhao~\cite{Wang13} found that
$h_{D_n}(-y^2)=\left(1+y^2\right)^n g_n(\phi)$, where
\begin{equation*}\label{gn}
  g_n(\phi)
=
\cos n\phi+
\frac n2\sin^2\phi\cos^{n-2}\phi.
\end{equation*}
Combining this with~\eqref{hD}, one can immediately get
\begin{equation*}
h_{C_n}(-y^2)=\left(1+y^2\right)^n\cos n\phi.
\end{equation*}
In order to establish the relationship among the zeros of $h_{C_n}(x)$ and $h_{D_n}(x)$, it suffices to explore that of 
$\cos n\phi$ and $g_n(\phi)$.

Following~\cite{Wang13}, by the inequality of arithmetic and geometric means, we obtain
\begin{align*}
\left(\frac{n}{2}\sin^2\phi\right)^2\left(\cos^{2}\phi\right)^{n-2} \leqslant\left(\frac{n\sin^2\phi+(n-2)\cos^{2}\phi}{n}\right)^n<1.
\end{align*}

Let
$
\alpha_i=\frac{(2i-1)\pi}{2n}$,
$ 1\leq i\leq n,
$
and
$
\beta_j=\frac{j\pi}{n}$,
$ 1\leq j\leq n-1.
$
Thus
$
\alpha_j<\beta_j<\alpha_{j+1},
$
and $\alpha_1,\ldots,\alpha_n$ are the zeros of
$\cos n\phi$. We consider the sign of $g_n(\alpha_j)$ and $g_n(\beta_j)$.
Note that
\begin{equation}\label{sign-alpha}
g_n(\alpha_j)
=
\frac n2\sin^2\alpha_j\cos^{n-2}\alpha_j
\end{equation}
and
\begin{equation}\label{sign-beta}
\operatorname{sgn}g_n(\beta_j)
=
\operatorname{sgn}\cos(n\beta_j)
=
(-1)^j.
\end{equation}
We now distinguish the parity of $n$.

\medskip
\noindent
\textbf{Case 1: $n=2m$.}

In this case, $n-2$ is even, and by \eqref{sign-alpha}, we have
$
g_n(\alpha_j)>0,
\ 1\leq j\leq n
$.
On the other hand, by \eqref{sign-beta}, when $j$ is odd, we have
$g_n(\beta_j)<0$.
Hence, for every odd $j$, by the intermediate value theorem, there exists
at least one zero of $g_n$ in each of the intervals
$(\alpha_j,\beta_j)$ and $(\beta_j,\alpha_{j+1})$.
For $j=1,3,\ldots,n-1$, we find $n$ zeros of $g_n$ in
$(0,\pi)$, which are all the zeros of $g_n$. Let $0<\theta_1<\theta_2<\cdots<\theta_n<\pi$ be the zeros of $g_n$.
Then
\begin{equation}\label{phi-even-order}
\alpha_1<\theta_1<\theta_2<\alpha_2
<\alpha_3<\theta_3<\theta_4<\alpha_4
<\cdots
<\alpha_{n-1}<\theta_{n-1}<\theta_n<\alpha_n.
\end{equation}

\medskip
\noindent
\textbf{Case 2: $n=2m+1$.}

Now $n-2$ is odd. It follows from \eqref{sign-alpha} that
$g_n(\alpha_j)>0,\  1\leq j\leq m$.
Moreover, $g_n\left(\frac{\pi}{2}\right)=0$.
Thus $\phi=\frac{\pi}{2}$ is the common zero of 
$\cos n\phi$ and $g_n(\phi)$, and in this case $x=-1$.
Furthermore,
\begin{equation}\label{symm}
  g_n(\pi-\phi)=\cos n(\pi -\phi) + \frac{n}{2} \sin^2(\pi -\phi) \cos^{n-2}(\pi -\phi) 
=-g_n(\phi),
\end{equation}
so the remaining zeros of $g_n$ occur symmetrically with
respect to $\frac{\pi}{2}$. Therefore it suffices to locate the zeros
in $(0,\frac{\pi}{2})$.

If $m$ is even, then for every odd
$j=1,3,\ldots,m-1$, we have
$g_n(\alpha_j)>0$,
$g_n(\beta_j)<0$, and
$g_n(\alpha_{j+1})>0$.
Hence there exists one zero in each of
$(\alpha_j,\beta_j)$ and $(\beta_j,\alpha_{j+1})$.
Thus we obtain $m$ zeros in $(0,\frac{\pi}{2})$, ordered as
\begin{equation}\label{phi-m-even}
\alpha_1<\theta_1<\theta_2<\alpha_2
<\alpha_3<\theta_3<\cdots
<\alpha_{m-1}<\theta_{m-1}<\theta_m<\alpha_m
<\alpha_{m+1}=\theta_{m+1}.
\end{equation}

If $m$ is odd, then the same argument applied to
$j=1,3,\ldots,m-2$ gives $m-1$ zeros in $(0,\frac{\pi}{2})$.
Moreover,
$g_n(\alpha_m)>0$ and $g_n(\beta_m)=(-1)^m=-1$.
Hence there exists an additional zero
$\theta_m\in(\alpha_m,\beta_m)$.
Therefore, the zeros in $(0,\frac{\pi}{2}]$ are ordered as
\begin{equation}\label{phi-m-odd}
\alpha_1<\theta_1<\theta_2<\alpha_2
<\alpha_3<\theta_3<\cdots
<\theta_{m-1}
<\alpha_{m-1}<\alpha_m<\theta_m
<\alpha_{m+1}=\theta_{m+1}.
\end{equation}

It follows from \eqref{symm} that $\theta_{n+1-j}=\pi-\theta_j$. Hence the inequalities~\eqref{phi-m-even} and \eqref{phi-m-odd} determine the order of all the zeros in $(0,\pi)$.
Finally, since the transformation $x=-\tan^2\frac{\phi}{2}$
is strictly decreasing on $(0,\pi)$, the zeros of $h_{C_n}(x)$ and
$h_{D_n}(x)$ are given by
\[
x_{n,j}^{c}
=
-\tan^2\frac{\alpha_{n+1-j}}{2},
\qquad
x_{n,j}^{d}
=
-\tan^2\frac{\theta_{n+1-j}}{2},
\qquad 1\leq j\leq n,
\]
and reversing the orderings
\eqref{phi-even-order}, \eqref{phi-m-even}, and
\eqref{phi-m-odd} yields
\eqref{interlace2m}, \eqref{interlace4r+1}, and
\eqref{interlace4r+3}, respectively.
This completes the proof.

\begin{example}
We have
\[
\begin{aligned}
h_{C_6}(x)
&=x^6+66x^5+495x^4+924x^3+495x^2+66x+1,\\
h_{D_6}(x)
&=x^6+54x^5+447x^4+852x^3+447x^2+54x+1.
\end{aligned}
\]
The zeros of $h_{C_6}(x)$ are approximately
\[
\begin{aligned}
-57.695481,\,-5.828427,\,-1.698396,
-0.588791,\,-0.171573,\,-0.017332,
\end{aligned}
\]
and the zeros of $h_{D_6}(x)$ are approximately
\[
\begin{aligned}
-44.348909,\,-7.207255,\,-1.691265,
-0.591273,\,-0.138749,\,-0.022548.
\end{aligned}
\]
Thus,
\[
\begin{aligned}
x_{6,1}^{c}
&<x_{6,1}^{d}
<x_{6,2}^{d}
<x_{6,2}^{c}
<x_{6,3}^{c}
<x_{6,3}^{d}
<x_{6,4}^{d}
<x_{6,4}^{c}
<x_{6,5}^{c}
<x_{6,5}^{d}
<x_{6,6}^{d}
<x_{6,6}^{c},
\end{aligned}
\]
which agrees with the interlacing relation in the even case.
\end{example}

\begin{example}
We have
\[
\begin{aligned}
h_{C_5}(x)
&=x^5+45x^4+210x^3+210x^2+45x+1,\\
h_{D_5}(x)
&=x^5+35x^4+180x^3+180x^2+35x+1.
\end{aligned}
\]
The zeros of $h_{C_5}(x)$ are approximately
\[
\begin{aligned}
-39.863458,\,-3.851840,\,-1,
-0.259616,\,-0.025086,
\end{aligned}
\]
and those of $h_{D_5}(x)$ are
\[
\begin{aligned}
-29.007120,\,-4.747781,\,-1,
-0.210625,\,-0.034474.
\end{aligned}
\]
Hence,
\[
\begin{aligned}
x_{5,1}^{c}
<x_{5,1}^{d}
<x_{5,2}^{d}
<x_{5,2}^{c}
<x_{5,3}^{c}
=x_{5,3}^{d}
=-1
<x_{5,4}^{c}
<x_{5,4}^{d}
<x_{5,5}^{d}
<x_{5,5}^{c},
\end{aligned}
\]
which agrees with the interlacing relation in the case
$n=2m+1$ with even $m$.
\end{example}

\begin{example}
We have
\[
\begin{aligned}
h_{C_7}(x)
&=x^7+91x^6+1001x^5+3003x^4
+3003x^3+1001x^2+91x+1,\\
h_{D_7}(x)
&=x^7+77x^6+931x^5+2863x^4
+2863x^3+931x^2+77x+1.
\end{aligned}
\]
The zeros of $h_{C_7}(x)$ are approximately
\[
\begin{aligned}
-78.769982,\,-8.167226,\,-2.532843,\,-1,
-0.394813,\,-0.122441,\,-0.012695,
\end{aligned}
\]
and the zeros of $h_{D_7}(x)$ are approximately
\[
\begin{aligned}
-62.913871,\,-10.071006,\,-2.499921,\,-1,
-0.400013,\,-0.099295,\,-0.015895.
\end{aligned}
\]
Therefore,
\[
\begin{aligned}
&x_{7,1}^{c}
<x_{7,1}^{d}
<x_{7,2}^{d}
<x_{7,2}^{c}
<x_{7,3}^{c}
<x_{7,3}^{d}
<x_{7,4}^{c}
=x_{7,4}^{d}
=-1\\
&<x_{7,5}^{d}
<x_{7,5}^{c}
<x_{7,6}^{c}
<x_{7,6}^{d}
<x_{7,7}^{d}
<x_{7,7}^{c},
\end{aligned}
\]
which agrees with the interlacing relation in the case
$n=2m+1$ with odd $m$.
\end{example}
%%%%%%%%%%%%%%%%%%%%%%%%%%%%%%%%%%%%%%%%%%%%%%%%%%%%%%%%%%%%%%%%%%%%%%%%%%
%%%%%%%%%%%%%%%%%%%%%%%%%%%%%%%%%%%%%%%%%%%%
\section*{Acknowledgements}
%\hspace*{\parindent}
%%%%%%%%%%%%%%%%%%%%%%%%%%%%%%%%%%%%%%%%%%%
This work is supported by the National Natural Science Foundation of China (No. 12071063) and 
the Natural Science Foundation of Shandong Province of China (ZR2026MS0047).
%%%%%%%%%%%%%%%%%%%%%%%%%%%%%%%%%%%%%%%%%%%%%%%%%%%%%%%%%%%%%%%%%%%%%%%%%%%%%%%%%
%%%%%%%%%%%%%%%%%%%%%%%%%%%%%%%%%%%%%%%%%%%%%%%%%

\end{document}